\documentclass[graybox]{svmult}
\usepackage{mathptmx}       
\usepackage{helvet}         
\usepackage{courier}        
\usepackage{makeidx}         
\usepackage{graphicx}
\usepackage{float}
\usepackage{xcolor}
\usepackage{booktabs} 
\usepackage{caption}
\usepackage{multicol}        
\usepackage[bottom]{footmisc}
\usepackage{amsmath,amssymb,amsfonts}        
\usepackage{bm}
\usepackage[normalem]{ulem} 

\makeindex             

\def\bds{\begin{displaystyle}}
\def\eds{\end{displaystyle}}
\def\E{\mathbb{E}}

\begin{document}
\title*{Estimation of rates in population-age-dependent  processes by means of test functions}
\author{Jie Yen Fan, Kais Hamza,  Fima C. Klebaner and Ziwen Zhong}
\institute{Jie Yen Fan \at School of Mathematics, Monash University, VIC 3800, Australia,  \\\email{Jieyen.Fan@monash.edu
\and Kais Hamza \at School of Mathematics, Monash University, VIC 3800, Australia, \\\email{Kais.Hamza@monash.edu}}
\and Fima C. Klebaner \at School of Mathematics, Monash University, VIC 3800, Australia, \\\email{Fima.Klebaner@monash.edu}
\and Ziwen Zhong \at School of Mathematics, Monash University, VIC 3800, Australia, \\\email{Ziwen.Zhong@monash.edu}}
%
%
\maketitle

\abstract{This paper aims to develop practical applications of the model for the highly technical measure-valued populations developed by the authors in \cite{FanEtal20}. We consider the problem of estimation of parameters in the general age and population-dependent model, in which the individual birth and death rates depend not only on the age of the individual but also on the whole population composition. We derive new estimators of the rates based on the use of test functions in the functional Law of Large Numbers and Central Limit Theorem for populations with a large carrying capacity.  
We consider the rates to be simple functions, that take finitely many values both in age $x$ and measure $A$, which leads to systems of linear equations. 
The proposed method of using time-dependent test functions for estimation is a novel approach which can be applied to a wide range of models of dynamical systems.
}


\section{Introduction}

In a sequence of papers \cite{JagKle00, HamJagKle16, FanEtal20, FanEtal20+} the authors developed a theory for general populations in which rates depend on the composition of the population as well as on the individual's age. This presents an important development in stochastic population dynamics theory. The evolution of population is determined by the way the individuals enter and the way they exit, which in turn are governed by the birth rate and the death rate respectively. It is mathematically convenient to describe the population as an atomic  measure $A$ on the line, and evolution in time as a measure-valued process $A_t$. These parameters $h$ and $b$ are assumed to depend on the age of the individual $x$ as well as on the population composition $A$.

The mathematically simplifying assumption is the introduction of the carrying capacity $K$,  which allows for approximations for large values of $K$. 
The results   in the above mentioned papers derive approximations for the composition of the population, which is  intractable otherwise.
The first approximation is the generalised McKendrick-von Foerster PDE, and the second approximation is for the fluctuations around it, given by a stochastic PDE (SPDE). 
We use these results for the estimation of unknown rates.

Since the analysis of measure-valued processes is done by means of test functions, it is natural to use them for estimation as well. 
The use of time-dependent test functions for estimation is a novel approach and yields new consistent estimators with a proven degree of accuracy. 
Moreover, in many cases, the problem reduces to solving a system of {\it linear} equations.

The numerical experiments back up our theory and show that this approach works. In particular, in the classical case of constant rates we recover the classical estimators \cite{Keiding74}.
 
Our work fits at the boundary between statistical 
learning and dynamical systems, in which parameters  are estimated
from the observed trajectory of dynamics equations. 

This work is the first step in developing inference by using time-dependent test functions, and has wide applicability in other areas.
Another advantage of our approach is the ability to estimate a multitude of parameters, by taking as many  test functions as necessary. 
This contrasts with the inability to estimate separately birth and death parameters in classical approach of birth-death process \cite{Keiding90} and particle kernel estimators \cite{Boumezoued21}, overcoming the problem of identifiability.

Population modelling and the estimation/ recovery of rates lie  at the intersection of many  areas.
Firstly, demography, where these rates are determined from mortality tables. 
Secondly, data driven models, in which various methods including statistics, are used to describe and explain observed population data.
Thirdly, population dynamics, where  models are based on the McKendrick-von Foerster PDE. 
Fourthly,  statistics, where underlying probability models are used, and is our approach.

Each of these areas has huge literature, and  here we mention just a few. Demographic models, see \cite{PreHeuGui01},  can be classified as data-driven models.  Population dynamics models based on analysis of the McKendrick-von Foerster PDE, eg. \cite{KeyKey97, Inaba17}  include the question of identifiability, ie. the ability to recover rates from observations, \cite{LiEtal21, Losanova22, RenKirEis22}. In statistical approach often populations are modelled by birth-death processes and branching processes. Their inference developed in \cite{Keiding74, Keiding74+, Keiding90, Guttorp91}.
   There are also studies on the age-dependent models, eg. \cite{OlaEtal24}.   The models in which  rates depend on age and population composition generalise branching models but technically they are not branching processes because the branching property is lost. Such models are also close to interacting particle systems.   A  model similar to ours is considered  in \cite{Boumezoued21},  where  kernel estimators are used to estimate the density function of the age process. However, none of these models consider rates that depend on both the age as well as population structure.

This work is the first step, as we mentioned already, and many 
questions remain, such as, which test functions to use, balancing mathematical and computational tractability on the one hand, with optimality, such as variance minimising, on the other.  We suggest that our approach  overcomes the problem of identifiability in general,  we demonstrate it for the model considered here, with more general statement to be addressed in further research.
We note that our estimators are consistent, due to the asymptotic theory developed earlier, however asymptotic normality is still to be established.
 
Section \ref{S:Prelim} formulates the general model, including the results on the Law of Large Numbers (LLN) and Central Limit Theorem (CLT). 
Section \ref{S:Equations} demonstrates how estimators of the parameters can be obtained from the LLN.
Numerical examples for some specific cases are also given. 
Section \ref{S:CI} explores the confidence intervals of the parameters using an auxiliary result of the CLT.

\section{Preliminaries} \label{S:Prelim}

We consider evolution of a population of finitely many individuals, whose ages we consider as a counting measure $A_t$ at time $t$ on $\mathbb{R}^+$,
$A_t(B)=\sum_{x\in A_t} 1_B(x)$. Here with a slight abuse of notation, we mean that $x$ is an atom of $A_t$ and $B$ is an interval. Each individual dies with   rate $h$ and gives birth with rate $b$. These parameters are assumed to depend on the age of the individual $x$ as well as on the population composition $A$, so that $h=h_A(x)$ and $b=b_A(x)$. Conditioned on the population composition, individuals act independently.
Furthermore, we assume large carrying capacity $K$, so that all the quantities are also indexed by $K$.
Our theory applies to populations evolving in time $t\in[0,T]$ for some arbitrary large but finite $T$.

It turns out that the measure-valued process $A^K_t$ is a Markov process with generator given in \cite{JagKle00}, from which \eqref{Evo1} can be derived.

For a $C^1$ function $f$ and a measure $A$, let $(f,A)=\int f(x)A(dx)$. Then the evolution equation is given by
\begin{equation}\label{Evo1}
(f,A^K_t)=(f,A_0^K)+\int_0^t (L^K_{A_s^K}f,A_s^K)ds +M^{K,f}_t,
\end{equation}  
where 
\begin{equation}\label{OperL}
L^K_Af=f'-h^K_Af+f(0)b^K_A
\end{equation}
is a first order differential operator and $M$ is a martingale. 
This equation was  generalised  for test functions that depend also on time, $f(x,t) \in C^{1,1}$,   \cite[Proposition 4]{FanEtal20}.
 
Writing $f_t(x)$ for $f(x,t)$,  
we have for any $t$ and $f\in C^{1,1}$
\begin{equation*} 
(f_t,A^K_t) = (f_0,A^K_0) + \int_0^t \big( L^K_{A^K_s}f_s ,A^K_s \big) ds
+ M^{K,f}_t,
\end{equation*}
where
\begin{equation}\label{OperL2}
L^K_{A}f(x,s) = \partial_1f(x,s) + \partial_2f(x,s) - f(x,s)h^K_{A} + f(0,s)b^K_{A}
\end{equation}
and $M_t^{K,f}$ is a  martingale with a known formula for its predictable quadratic variation. 

While we use the same notation $L^K_A$ for the operator in both equations, it is clear from the context which of \eqref{OperL} or \eqref{OperL2} applies.

We further assume that as $K\to\infty$ the  parameters ($b^K$ and $h^K$) tend to their limiting values, functions (of population $A$ and age $x$) $h_A(x)$ and $b_A(x)$, forming conditions  we termed {\it  smooth demography} in \cite{FanEtal20}. 
This paper aims to estimate $h_A(x)$ and $b_A(x)$.

It is shown in \cite{HamJagKle16}, see also \cite{FanEtal20}, that in smooth demographics  the functional LLN holds,  $\bar A^K_t := \frac{1}{K}A_t^K$ converges weakly to $\bar A_t$ (in appropriate Skorohod space of trajectories with values in space of positive measures). 
The limit process $\bar A_t$ is a deterministic measure satisfying equation
\begin{equation}\label{EvoA_2var}
(f_t,\bar A_t) = (f_0,\bar A_0) + \int_0^t (\partial_x f_s + \partial_t f_s - f_s h_{\bar A_s} + f_s(0) b_{\bar A_s}, \bar A_s) ds,
\end{equation}
where $\bar A_0$ is the limit  as $K\to\infty$ of  $\bar A^K_0 := \frac{1}{K}A_0^K$.
In particular, taking $f$ as a function of the first variable $x$ only, we have
\begin{equation}\label{EvoA}
(f,\bar A_t) = (f,\bar A_0) + \int_0^t (f' - f h_{\bar A_s} + f(0) b_{\bar A_s}, \bar A_s) ds.
\end{equation} 

Note that for practical applications of the model one can take $K$ as the size of  the initial population.

One can view  \eqref{EvoA_2var}  and  \eqref{EvoA}  as a weak form of the generalised McKendrick-von Foerster PDE.
The density $a(x,t)$ of $\bar A_t$ (with respect to Lebesgue measure)  solves the familiar McKendrick-von Foerster PDE \cite{Hoppenstaedt75,  McKendrick25, vonFoerster59}, but now it is generalised in the sense of allowing parameters $h$ and $b$ to  depend also on $A$, which makes the PDE into a non-linear one:
\begin{equation*}
\left(\frac{\partial}{\partial x}+\frac{\partial}{\partial t}\right)a(x,t) = -a(x,t)h_{\bar A_t}(x),\quad a(0,t)=\int_0^\infty b_{\bar A_t}(x) a(x,t) dx.
\end{equation*}

To obtain confidence bounds for the parameters, we use functional CLT for measure-valued populations obtained in \cite{FanEtal20}.
Under appropriate broad assumptions, the fluctuation process $Z_t^K := \sqrt K(\bar A_t^K- \bar A_t)$ converges (in the appropriate Skorohod space of trajectories with values in Sobolev space $W^{-4}$)  to a limit $Z$ satisfying an SPDE, \cite{FanEtal20}. We shall use an auxiliary fact in the proof of CLT \cite[Proposition 26]{FanEtal20} that the martingales $M^{K,f}_t$ in \eqref{Evo1} scaled by $\frac{1}{\sqrt K}$  converge to the Gaussian martingale $M^f_t$ with zero mean and quadratic variation
\begin{equation}\label{QVM}
\langle M^f,M^f \rangle_t=\int_0^t (f^2(0)b_{\bar A_s}+h_{\bar A_s}f^2, \bar A_s)ds.
\end{equation}

\section{Estimating Equations} \label{S:Equations}

The idea is to use the limiting evolution equation with various test functions to extract information about the rates.
Rearranging equation \eqref{EvoA} for parameters we obtain the starting point for their inference.

We have 
\begin{equation}\label{Main}
f(0)\int_0^t(b_{\bar A_s}, \bar A_s)ds - \int_0^t (h_{\bar A_s}f, \bar A_s)ds = (f,\bar A_t) - (f,\bar A_0)-\int_0^t(f',\bar A_s)ds.
\end{equation}
In some cases, we need a richer class of test functions, test functions that depend also on time, $f(x,t)$, written as $f_t(x)$ below. 
Equation \eqref{EvoA_2var}  gives
\begin{equation}\label{Main2}
\int_0^t f_s(0)(b_{\bar A_s}, \bar A_s)ds - \int_0^t (h_{\bar A_s}f_s, \bar A_s)ds = (f_t,\bar A_t) - (f_0,\bar A_0) - \int_0^t (\partial_x f_s + \partial_t f_s,\bar A_s)ds.
\end{equation}

Of course, if the limit $\bar A_t$ is known, no estimation is required as we can recover rates  by solving the above equations (inverse problem). We assume, however, that we observe the  pre-limit process $\bar A_t^K$, $0\le t\le T$, and  
the we obtain estimators by replacing the limit process $\bar A_t$ by its pre-limit $\bar A_t^K$ for large $K$. 
Note that this approach produces consistent  estimators  (as $K\to \infty$),  which follows from the weak  convergence of $\bar A_t^K$ to $\bar A$ (given by our LLN), and the Slutzky theorem  (which states that convergence is preserved under continuous transformations).

From \eqref{Main} or \eqref{Main2}, taking test functions that are null at 0 eliminates $b$ from the equation, leaving only $h$. This allows one to obtain $h$ first, and then obtain $b$. 

In what follows we consider models with increasing complexity, starting with constant parameters and ending with parameters fully dependent on the population as well as age of the individual. We consider  the rates to be simple  functions of its variables taking finitely many values both in age $x$ and measure $A$. This assumption   leads to systems of linear equations for  recovery of the constants. To justify this choice, note that from theoretical perspective, simple functions approximate any measurable function; and  
from practical perspective, it is intuitively clear that  one can assume the rates to be constants on various age intervals. In this first work on our new approach we don’t discuss how to choose the intervals of constancy, and leave this choice for later research. Our numerical examples provide conceptual training models to demonstrate the effectiveness and main steps of the approach.
Having said this, our approach is clearly applicable to other models of rates.

Regarding notations, we agree to write, with a slight abuse of notation, $(x,A)$ instead of  $(f,A)$ when $f(x)=x$, and  $(xt,A)$ when $f(x,t)=xt$, similarly for other explicit forms of test functions $f$.

\subsection{Constant parameters} \label{S:const}

Consider first the classical case of constant parameters $h$ and $b$, constant both in $x$ and $A$.
Then equation \eqref{Main} yields
\begin{equation*}
f(0) b \int_0^T (1, \bar A_s) ds - h\int_0^T (f, \bar A_s) ds = (f,\bar A_T) - (f,\bar A_0) - \int_0^T (f',\bar A_s) ds.
\end{equation*}
We take $f(x)=x$ to obtain $h$, (with $f(0)=0$, $f'(x)=1$)
\begin{equation*}
h = \frac{(x, \bar A_0) - (x,\bar A_T) + \int_0^T (1,\bar A_s) ds}{\int_0^T (x, \bar A_s) ds}.
\end{equation*}
Taking $f(x)=1$ we then obtain $b$, (with $f(0)=1$, $f'(x)=0$)
\begin{equation*}
b = \frac{(1,\bar A_T) - (1,\bar A_0) + h \int_0^T (1,\bar A_s) ds}{\int_0^T (1,\bar A_s) ds}.
\end{equation*}
Replacing the limit process $\bar A$ with $\bar A^K$, we obtain the estimators of $h$ and $b$.
 We simulate the full sample path of the measure-valued process. This allows us to evaluate the integrals in the estimating equations accurately, as these integrals depend on the continuous-time trajectory of the process.

\begin{figure}[h]
    \centering
    \includegraphics[width=9cm]{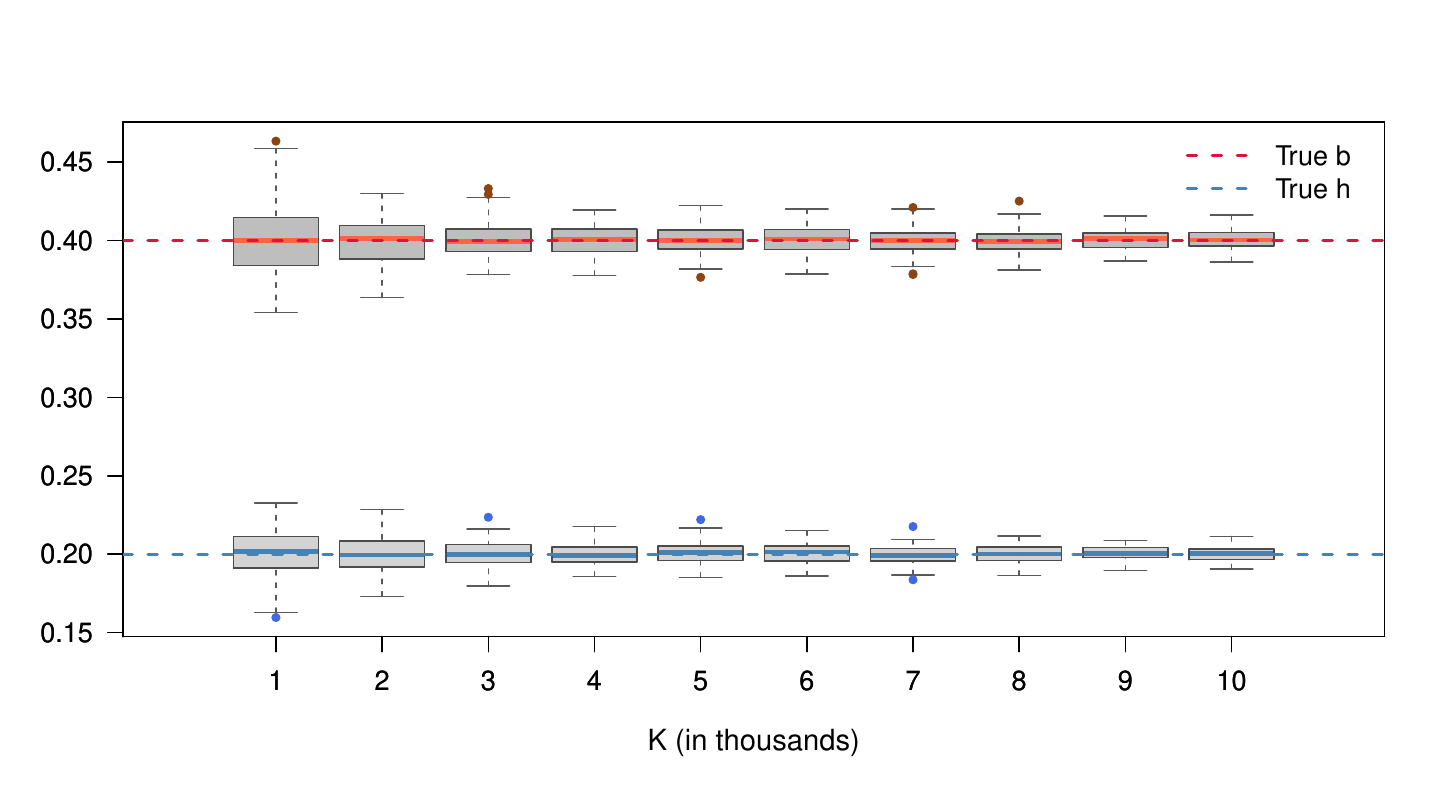}
    \caption{Box plots of 100 estimates of $b$ and $h$ with different $K$.} 
    \label{fig:constant}
\end{figure}

Numerical results are presented below with parameters \( h = 0.2 \), \( b = 0.4 \), and different initial population sizes $K$. 
The age of each individual at time 0 is taken to be randomly distributed in the interval $[0,1]$ following the uniform distribution, and $T=1$.
Tables \ref{tab:h_summary} and \ref{tab:b_summary} show some summary statistics of 100 estimates of $h$ and $b$ for different $K$.
Figure \ref{fig:constant} displays box plots of 100 estimates of $h$ and $b$ for different $K$. 

\begin{table}[h]
    \centering
    \begin{minipage}{0.45\textwidth}
        \centering
        \begin{tabular}{lccc}
            \toprule
            K & 100 & 1000 & 10000 \\
            \midrule
            Sample Mean    & 0.20874  & 0.19843  & 0.19962  \\
            Sample Variance & 0.00233 & 0.00023 & 0.00003 \\
            MSE            & 0.00238 & 0.00023 & 0.000027 \\
            Bias           & 0.00874 & -0.00157 & -0.00038 \\
            \bottomrule
        \end{tabular}
        \caption{Summary statistics of 100 estimates of $h$ with different $K$.}
        \label{tab:h_summary}
    \end{minipage}
    \hfill
    \begin{minipage}{0.45\textwidth}
        \centering
        \begin{tabular}{lccc}
            \toprule
            K & 100 & 1000 & 10000 \\
            \midrule
            Sample Mean    & 0.39848  & 0.39848  & 0.39878  \\
            Sample Variance & 0.00345 & 0.00040 & 0.00004 \\
            MSE            & 0.00342 & 0.00040 & 0.00004 \\
            Bias           & -0.00152 & -0.00152 & -0.00122 \\
            \bottomrule
        \end{tabular}
        \caption{Summary statistics of 100 estimates of $b$ with different $K$.}
        \label{tab:b_summary}
    \end{minipage}
\end{table}

\subsection{Parameters depend only on population} \label{S:PopDep}

Consider next the case where parameters $h$ and $b$ are constant in $x$ but depend on $A$. 
Equation \eqref{Main} yields in this case
\begin{equation*}
f(0) \int_0^T b_{\bar A_s}(1,\bar A_s) ds - \int_0^T h_{\bar A_s}(f,\bar A_s) ds = (f,\bar A_T) - (f,\bar A_0) - \int_0^T(f',\bar A_s) ds.
\end{equation*}

Some fairly general cases of functions of a measure include some function 
applied to the linear function of $A$, which
is $(\phi,A)$ for some $\phi$, $g((\phi,A))$, and such sums $\sum_{i,j} g_i((\phi_j,A))$.

However, for the purpose of modelling 
it is plausible that the dependence of the birth parameter on population is proportional to the number of individuals in a particular age interval $J_1$, i.e.  $b_A = \eta (1_{J_1},A)$ for some constant $\eta$. Similarly, the death parameter may be proportional to the number of individuals  in another  age interval $J_2$, i.e. $h_A = \lambda (1_{J_2},A)$ for some constant $\lambda$. 

Taking $f(x)=x$ and $f(x)=1$ allows us to obtain the following formulae for $\lambda$ and $\eta$:
\begin{equation*}
\lambda = \frac{(x, \bar A_0) - (x, \bar A_T) + \int_0^{T} (1,\bar A_s) ds}{\int_{0}^{T} (x,\bar A_s) (1_{J_2},\bar A_s) ds},
\end{equation*}
and
\begin{equation*}
\eta = \frac{(1,\bar A_T) - (1, \bar A_0) + \int_{0}^{T} \lambda (1_{J_2},\bar A_s)(1,\bar A_s)ds}{\int_{0}^{T}(1,\bar A_s)(1_{J_1},\bar A_s)ds}.
\end{equation*}
Replacing the limit process $\bar A$ with $\bar A^K$, we obtain the estimators of $\lambda$ and $\eta$.

\begin{figure}[H]
    \centering
    \includegraphics[width=9cm]{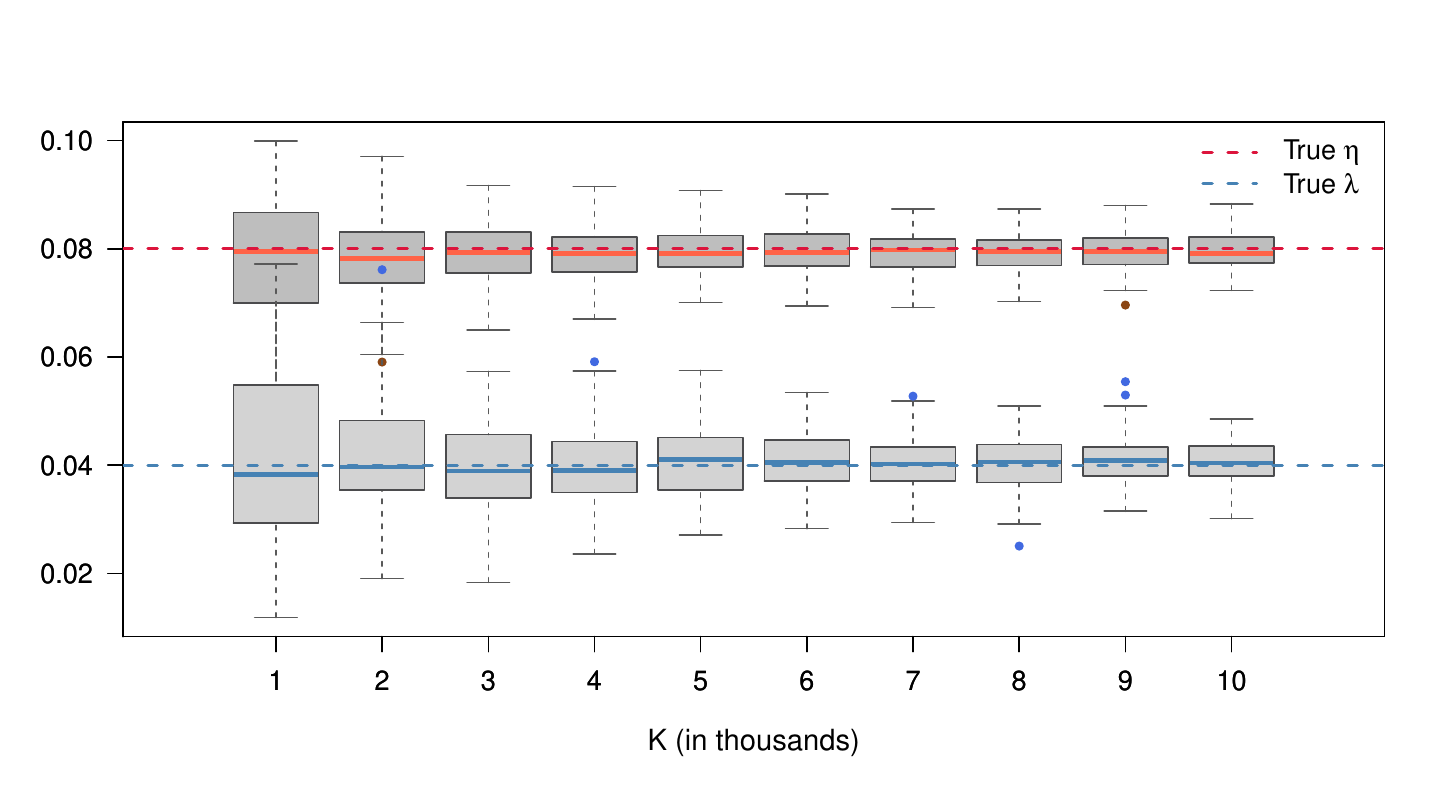}
    \caption{Box plots of 100 estimates of $\eta$ and $\lambda$ with different $K$.}
    \label{fig:PopDep}
\end{figure}

For example, take $J_1=[0.5, 1.5]$, $J_2= 1_{[0, 0.5) \cup (1.5, 2]}$, $\eta=0.08$, and $\lambda=0.04$, i.e. 
$$b_A = 0.08 ( 1_{[0.5, 1.5]},A) \quad \textnormal{and} \quad h_A = 0.04 (1_{[0, 0.5) \cup (1.5, 2]},A).$$
Let the age of each individual at time 0 follow the uniform distribution on $[0,1]$, and take $T=1$. 
We obtain the following numerical results from 100 sample paths for each chosen value of $K$.
Tables \ref{tab:d_summary} and \ref{tab:c_summary} show summary statistics of the 100 estimates of $\lambda$ and $\eta$ for different $K$.
Figure \ref{fig:PopDep} shows box plots of 100 estimates of $\lambda$ and $\eta$ for different $K$. 

\begin{table}[H]
    \centering
    \begin{minipage}{0.45\textwidth}
        \centering
        \begin{tabular}{lccc}
            \toprule
            K & 100 & 1000 & 10000 \\
            \midrule
            Sample Mean  & 0.04897  & 0.04174  & 0.04053  \\
            Sample Variance & 0.00181 & 0.00025 & 0.00001 \\
            MSE        & 0.00187 & 0.00025 & 0.00001 \\
            Bias          & 0.00897 & 0.00174 & 0.00053 \\
            \bottomrule
        \end{tabular}
        \caption{Summary statistics of 100 estimates of $\lambda$ with different $K$.}
        \label{tab:d_summary}
    \end{minipage}
    \hfill
    \begin{minipage}{0.45\textwidth}
        \centering
        \begin{tabular}{lccc}
            \toprule
            K & 100 & 1000 & 10000 \\
            \midrule
            Sample Mean    & 0.08479  & 0.07837  & 0.07956  \\
            Sample Variance & 0.00125 & 0.00015 & 0.00001 \\
            MSE            & 0.00126 & 0.00015 & 0.00001 \\
            Bias           & 0.00479 & -0.00163 & -0.00044 \\
            \bottomrule
        \end{tabular}
        \caption{Summary statistics of 100 estimates of $\eta$ with different $K$.}
        \label{tab:c_summary}
    \end{minipage}
\end{table}

Clearly, other explicit dependencies on $A$ can be incorporated in a similar way.

\subsection{Parameters depend only on age} \label{sec: age_depend}

Consider next the  case when parameters $h$ and $b$, are constant in $A$ but depend on $x$.
In this case we consider piecewise constant functions. While not the most general, bear in mind that such function approximate very wide class of functions of $x$.
It is natural to take 
\begin{equation*}
h(x)=\sum_{i=1}^n h_i1_{B_i}(x) \quad \textnormal{and} \quad b(x)=\sum_{i=1}^n b_i1_{B_i}(x), 
\end{equation*}
where $B_i$'s  are intervals (sets) on which parameters are constants. Of course, there is no essential difficulty to take different intervals of constancy for $b$ and $h$, but it seems make sense to say that $b$ and $h$ are constant in same age classes.

To recover $h$ we use \eqref{Main2} with functions $f_t(x)=xt^m$,  $m=0,1,2,\ldots,n-1$. 
Note that $f_t(0)=0$, $\partial_x f_t=t^m$, $\partial_t f_t=mxt^{m-1}$, and $f_0= x1_{m=0}$.
In this case, 
$$(h_{A}f_t, A) = \Big(\sum_{i=1}^n h_i1_{B_i}f_t,A\Big) = \sum_{i=1}^n h_i(1_{B_i}f_t,A).$$
Further, for $f_t(x) = xt^m$,
$(1_{B_i}f_t,A) = t^m(x1_{B_i}(x),A)$, giving
$$\int_0^T (h_{A_s}f_s, A_s) ds = \sum_{i=1}^n h_i\int_0^T s^m (x1_{B_i}(x),A_s) ds.$$
Thus we obtain a system of $n$ linear equations for $h_i$'s. 
For $m=0$,
\begin{equation}\label{AgeDepEq_h1}
\sum_{i=1}^n h_i \int_0^T (x1_{B_i}(x),\bar A_s) ds = (x, \bar A_0) - (x,\bar A_T) + \int_0^T (1,\bar A_s)ds,
\end{equation}
and for $m= 1, 2, \dots, n-1$,
\begin{equation}\label{AgeDepEq_h2}
\sum_{i=1}^n h_i \int_0^T s^m(x1_{B_i}(x),\bar A_s) ds = - T^m(x,\bar A_T) + \int_0^T s^m (1,\bar A_s) ds + m \int_0^T s^{m-1} (x,\bar A_s) ds.
\end{equation}
Denote $g_i(s)=(x1_{B_i}(x),\bar A_s)$. For any positive integer $n$, we write $[n]:=\{1, 2, \cdots, n\}$  to denote the set of the first 
$n$ natural numbers. The determinant of the matrix with elements 
$( \int_0^T s^m g_i(s)ds)_{i\in [n],\;m\in [n-1] \cup \{0\}}$  is not zero in general, which assures a unique solution.

Having found $h_i$'s we recover $b_i$'s next. To this end we use functions $f_t(x) = t^m$, 
i.e. $f_t(0) = t^m$,  $\partial_x f_t = 0$, $\partial_t f_t = mt^{m-1}$, and $f_0 = 1_{m=0}$. Note that 
$$(b_A,A) = \sum_{i=1}^n b_i(1_{B_i},A)$$
and
$$\int_0^T f_s(0) (b_{A_s}, A_s) ds = \sum_{i=1}^n b_i \int_0^T s^m (1_{B_i},A_s) ds.$$
Thus we obtain a system of $n$ linear equations for $b_i$'s. 
For $m=0$,
\begin{equation}\label{AgeDepEq_b1}
\sum_{i=1}^n b_i \int_0^T (1_{B_i},\bar A_s) ds = (1,\bar A_T) - (1, \bar A_0) + \sum_{i=1}^n h_i \int_0^T (1_{B_i},\bar A_s) ds,
\end{equation}
and for $m= 1, 2, \dots, n-1$,
\begin{equation}\label{AgeDepEq_b2}
\sum_{i=1}^n b_i\int_0^T s^m(1_{B_i},\bar A_s)ds = T^m(1,\bar A_T) - m \int_0^T s^{m-1} (1,\bar A_s)ds
+\sum_{i=1}^n h_i\int_0^T s^m(1_{B_i},\bar A_s)ds.
\end{equation}
Similar to the system for $h_i$'s, one can see that this system has a unique solution, by checking the non-degeneracy of its determinant.

The estimators of the parameters are then obtained by replacing the limit process $\bar A$ with $\bar A^K$.

\begin{figure}[h]
    \centering
    \begin{minipage}{0.43\textwidth} 
        \centering
        \includegraphics[width=6cm]{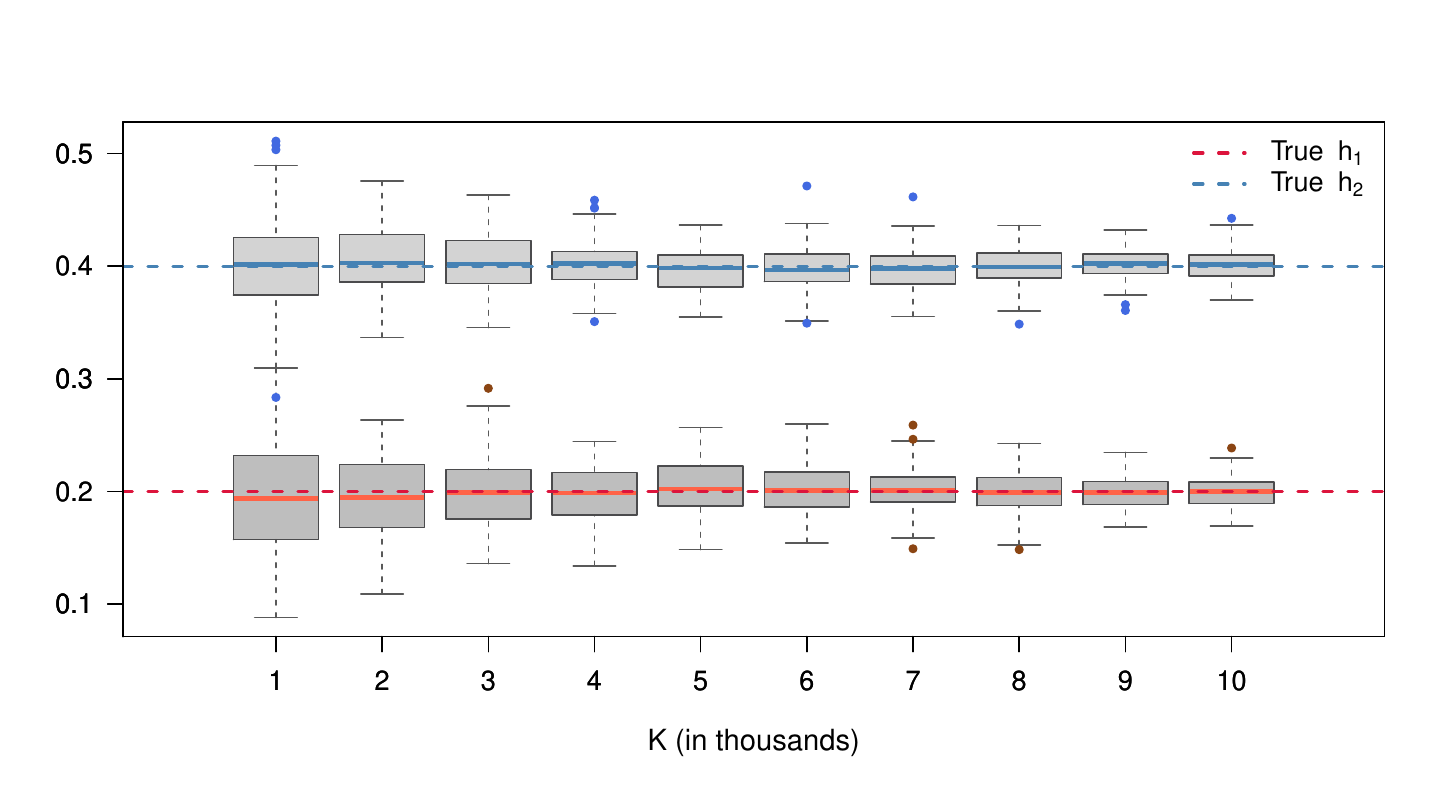} 
    \end{minipage}
    \hfill 
    \begin{minipage}{0.5\textwidth}
        \centering
        \includegraphics[width=6cm]{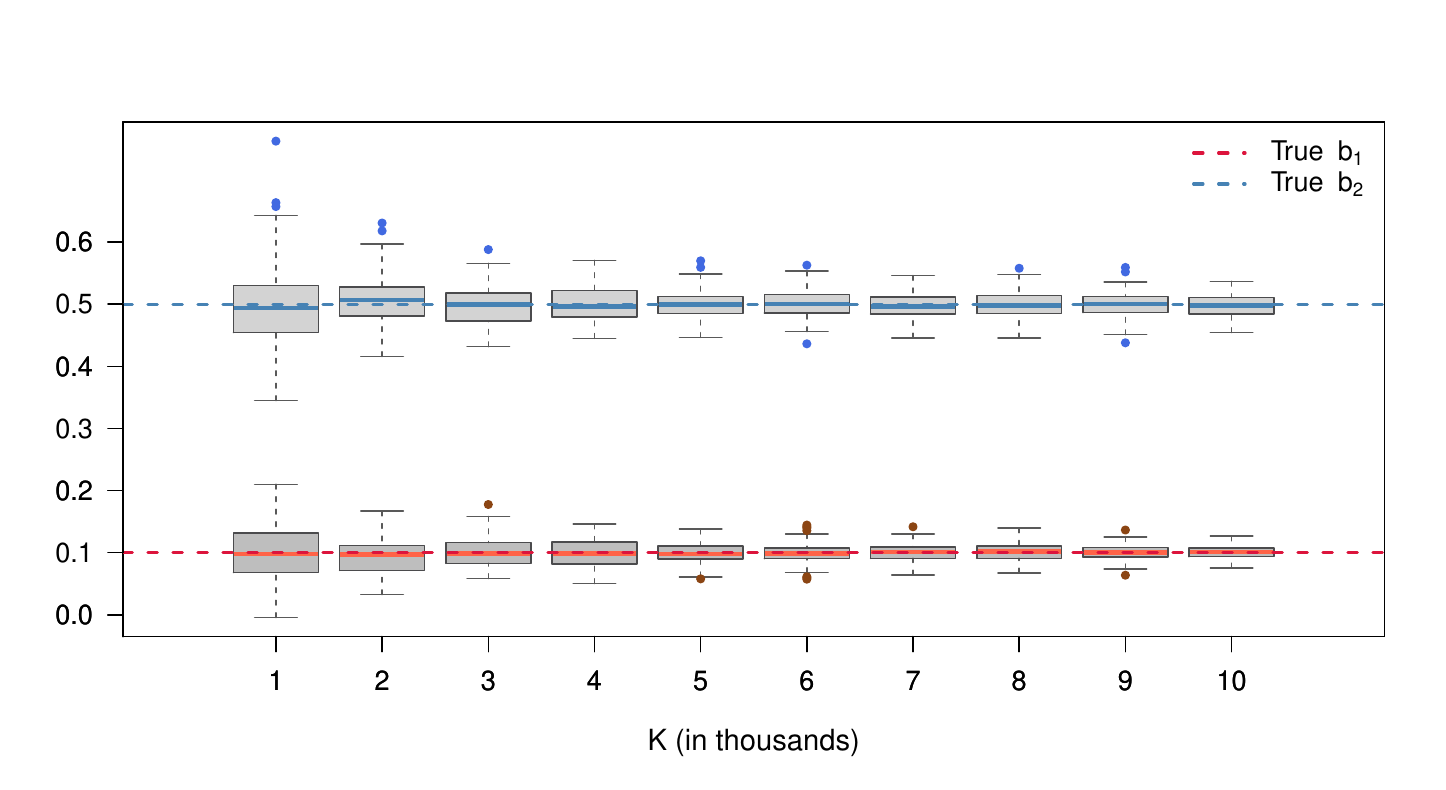} 
    \end{minipage}
    \caption{Box plots of 100 estimates of \( h_1 , h_2 \) (left) and \( b_1 , b_2 \) (right) with different $K$.}
    \label{fig:AgeDep_hb}
\end{figure}

Since in principle there is little difference between the case of $n=2$ and larger $n$ (except for computing time), we consider a numerical example for $n=2$.
Take $B_1=[0,1)$, $B_2=[1,2]$, $h_1=0.2$, $h_2=0.4$, $b_1=0.1$, and $b_2=0.5$, i.e. 
$$h(x) = 0.2 \ 1_{[0,1)}(x) + 0.4 \ 1_{[1,2]}(x) 
\quad \textnormal{and} \quad 
b(x) = 0.1 \ 1_{[0, 1)}(x) + 0.5 \ 1_{[1,2]}(x).$$
Suppose the age of each individual at time 0 is uniformly distributed on $[0,1]$. Take $T=1$.
With 100 sample paths for each chosen $K$ value, we obtain the following  results from Equations \eqref{AgeDepEq_h1}-\eqref{AgeDepEq_b2}. 
Tables \ref{tab:h_estimates} and \ref{tab:b_estimates} show summary statistics of 100 estimates for different $K$.
Figure \ref{fig:AgeDep_hb} shows box plots of 100 estimates of $h_1$, $h_2$, $b_1$, and $b_2$ for different $K$. 

\begin{table}[h]
    \centering
    \begin{tabular}{lcccccc}
        \toprule
        $K$ & \multicolumn{2}{c}{100} & \multicolumn{2}{c}{1000} & \multicolumn{2}{c}{10000} \\
        \cmidrule(lr){2-3} \cmidrule(lr){4-5} \cmidrule(lr){6-7}
        & $h_1$ & $h_2$ & $h_1$ & $h_2$ & $h_1$ & $h_2$ \\
        \midrule
        Sample Mean & 0.18009 & 0.40605 & 0.19771 & 0.40167 & 0.19991 & 0.40106 \\
        Sample Variance & 0.03038 & 0.02220 & 0.00270 & 0.00159 & 0.00016 & 0.00019 \\
        MSE & 0.03047 & 0.02201 & 0.00268 & 0.00158 & 0.00016 & 0.00019 \\
        Bias & -0.01913 & 0.00605 & -0.00229 & 0.00167 & -0.00009 & 0.00106 \\
        \bottomrule
    \end{tabular}
    \caption{Summary statistics of 100 estimates of $h_1$ and $h_2$ with different $K$.}
    \label{tab:h_estimates}
\end{table}

\begin{table}[h]
    \centering
    \begin{tabular}{lcccccc}
        \toprule
        $K$ & \multicolumn{2}{c}{100} & \multicolumn{2}{c}{1000} & \multicolumn{2}{c}{10000} \\
        \cmidrule(lr){2-3} \cmidrule(lr){4-5} \cmidrule(lr){6-7}
        & $b_1$ & $b_2$ & $b_1$ & $b_2$ & $b_1$ & $b_2$ \\
        \midrule
        Sample Mean & 0.08658 & 0.51625 & 0.09803 & 0.49722 & 0.10087 & 0.49802 \\
        Sample Variance & 0.01933 & 0.03466 & 0.00199 & 0.00436 & 0.00012 & 0.00033 \\
        MSE & 0.01932 & 0.03458 & 0.00197 & 0.00432 & 0.00012 & 0.00033 \\
        Bias & 0.01342 & 0.01625 & -0.00197 & -0.00278 & 0.00087 & -0.00198 \\
        \bottomrule
    \end{tabular}
    \caption{Summary statistics of 100 estimates of $b_1$ and $b_2$ with different $K$.}
    \label{tab:b_estimates}
\end{table}

\subsection{Parameters depend on population and age} \label{sec: popage}

Consider the general case when parameters $h$ and $b$ depend on both $A$ and $x$. Again, we consider a somewhat simplified situation
when dependence on $x$ is piecewise constant, i.e. 
\begin{equation*}
h_A(x) = \sum_{i=1}^n h^{(i)}_A 1_{B_i}(x) \quad \textnormal{and} \quad b_A(x) = \sum_{i=1}^n b^{(i)}_A 1_{B_i}(x),
\end{equation*}
where $h^{(i)}_A$ and $b^{(i)}_A$ are constant in $x$ but depend on $A$.
From equation \eqref{Main2},
\begin{multline*}
\sum_{i=1}^{n} \int_0^T b^{(i)}_{\bar A_s}f_s(0) (1_{B_i}, \bar A_s) ds - \sum_{i=1}^{n} \int_0^T h^{(i)}_{\bar A_s} (f_s 1_{B_i}, \bar A_s) ds \\
= (f_T,\bar A_T) - (f_0,\bar A_0) - \int_0^T (\partial_x f_s + \partial_t f_s, \bar A_s) ds.
\end{multline*}
Similar to the approach considered in Section \ref{sec: age_depend}, using test functions $f_t(x)=xt^m$ and $f_t(x)=t^m$, for $m=0,1,2,\ldots,n-1$, we can recover $h^{(i)}_{A}$ and $b^{(i)}_{A}$.

For example, let
\begin{equation} \label{eq: popage_h}
h_A(x)=\alpha_1 (1_{J},A) 1_{B_1}(x) + \alpha_2 (1_{J},A) 1_{B_2}(x),
\end{equation}
and
\begin{equation} \label{eq: popage_b}
b_A(x)=\gamma_1 (1_{J},A )1_{B_1}(x) + \gamma_2 (1_{J},A) 1_{B_2}(x). 
\end{equation}
Taking $f_t(x) = x$ and $f_t(x) = xt$, we can recover $\alpha_1$ and $\alpha_2$ by solving 
\begin{equation*}
\sum_{i=1}^2 \alpha_i \int_0^T (1_{J}, \bar A_s)(x 1_{B_i}(x),\bar A_s) ds = (x, \bar A_0) - (x,\bar A_T) + \int_0^T (1,\bar A_s)ds
\end{equation*}
and
\begin{equation*}
\sum_{i=1}^2 \alpha_i \int_0^T s (1_{J}, \bar A_s)(x1_{B_i}(x),\bar A_s)ds = - T(x,\bar A_T) + \int_0^T s(1,\bar A_s)ds + \int_0^T (x,\bar A_s)ds.
\end{equation*}
Having found $\alpha_i$'s, we can recover $\gamma_i$'s next. 
Taking $f_t(x) = 1$ and $f_t(x) = t$, we have 
\begin{equation*}
\sum_{i=1}^2 \gamma_i \int_0^T (1_{J},\bar A_s)(1_{B_i}(x), \bar A_s) ds 
= (1,\bar A_T) - (1, \bar A_0) + \sum_{i=1}^2 \alpha_i \int_0^T (1_J, \bar A_s)(1_{B_i}(x),\bar A_s) ds
\end{equation*}
and 
\begin{multline*}
\sum_{i=1}^2 \gamma_i\int_0^T s(1_J, \bar A_s)(1_{B_i}(x), \bar A_s)ds \\
= T(1,\bar A_T) - \int_0^T (1,\bar A_s) ds + \sum_{i=1}^n \alpha_i \int_0^T s (1_J, \bar A_s)(1_{B_i}(x), \bar A_s) ds.
\end{multline*}
Replacing the limit process $\bar A$ with $\bar A^K$ we obtain the estimators of $\alpha_i$'s and $\gamma_i$'s.

\begin{figure}[h]
    \centering
    \begin{minipage}{0.43\textwidth}
        \centering
        \includegraphics[width=6cm]{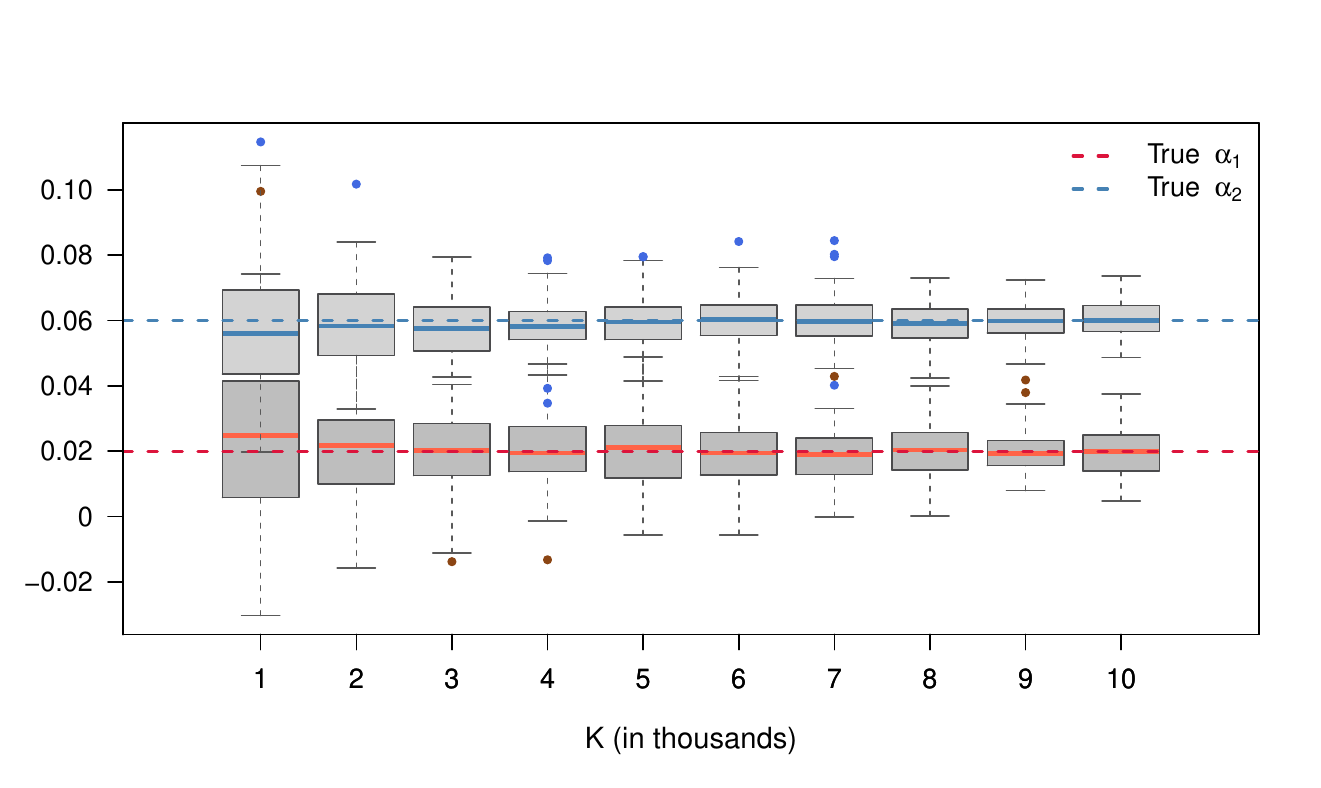} 
    \end{minipage}
    \hfill 
    \begin{minipage}{0.5\textwidth}
        \centering
        \includegraphics[width=6cm]{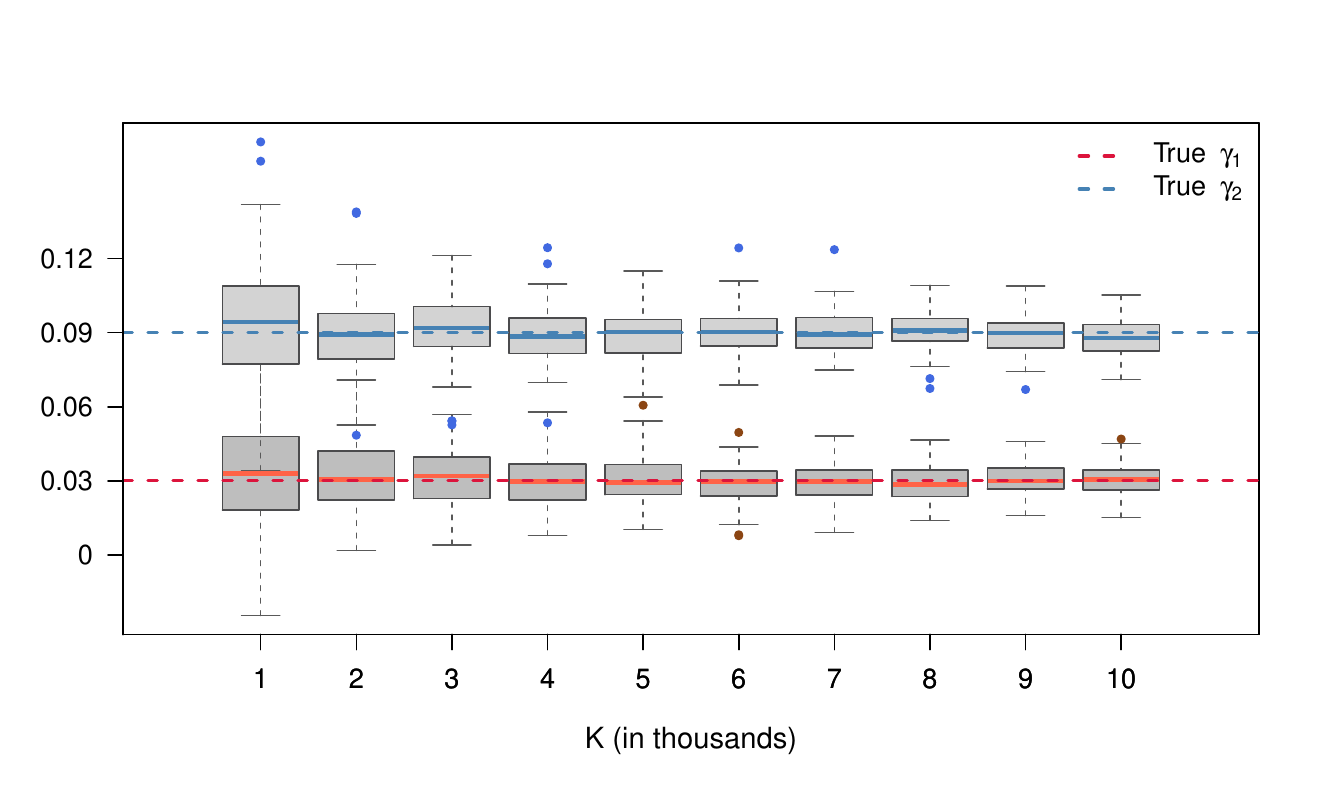} 
    \end{minipage}
    \caption{Box plots of 100 estimates of \( \alpha_1 , \alpha_2 \) (left) and \( \gamma_1 , \gamma_2 \) (right) with different $K$.}
    \label{fig:PopAgeDep}
\end{figure}

For numerical example, we take $J = [0.5, 1.5]$, $B_1 = [0,1)$, $B_2 = [1, 2]$, 
$\alpha_1 = 0.02$, $\alpha_2 = 0.06$, $\gamma_1= 0.03$, $\gamma_2 = 0.09$. 
As before, the age of each individual at time 0 is taken to follow uniform distribution on $[0,1]$, and $T=1$. 
With 100 sample paths for each chosen value of $K$, we obtain the following numerical results. 
Tables \ref{tab:alpha_estimates} and \ref{tab:beta_estimates} show summary statistics of 100 estimates of $\alpha_i$ and $\gamma_i$ for different $K$. 
Figure \ref{fig:PopAgeDep} shows box plots of 100 estimates of the $\alpha_i$'s and $\gamma_i$'s for different $K$. 

\begin{table}[h]
    \centering
    \begin{tabular}{lcccccc}
        \toprule
        $K$ & \multicolumn{2}{c}{100} & \multicolumn{2}{c}{1000} & \multicolumn{2}{c}{10000} \\
        \cmidrule(lr){2-3} \cmidrule(lr){4-5} \cmidrule(lr){6-7}
        & $\alpha_1$ & $\alpha_2$ & $\alpha_1$ & $\alpha_2$ & $\alpha_1$ & $\alpha_2$ \\
        \midrule
        Sample Mean & 0.02373 & 0.06684 & 0.02348 & 0.05778 & 0.01968 & 0.06029 \\
        Sample Variance & 0.00581 & 0.00389 & 0.00062 & 0.00038 & 0.00005 & 0.00003 \\
        MSE & 0.00577 & 0.00390 & 0.00063 & 0.00038 & 0.00005 & 0.00003 \\
        Bias & 0.00373 & 0.00684 & 0.00348 & -0.00222 & -0.00032 & 0.00029 \\
        \bottomrule
    \end{tabular}
    \caption{Summary statistics of 100 estimates of $\alpha_1$ and $\alpha_2$ with different $K$.}
    \label{tab:alpha_estimates}
\end{table}

\begin{table}[h]
    \centering
    \begin{tabular}{lcccccc}
        \toprule
        $K$ & \multicolumn{2}{c}{100} & \multicolumn{2}{c}{1000} & \multicolumn{2}{c}{10000} \\
        \cmidrule(lr){2-3} \cmidrule(lr){4-5} \cmidrule(lr){6-7}
        & $\gamma_1$ & $\gamma_2$ & $\gamma_1$ & $\gamma_2$ & $\gamma_1$ & $\gamma_2$ \\
        \midrule
        Sample Mean & 0.02677 & 0.09750 & 0.03221 & 0.09362 & 0.03053 & 0.08843 \\
        Sample Variance & 0.00349 & 0.00614 & 0.00047 & 0.00056 & 0.00004 & 0.00005 \\
        MSE & 0.00346 & 0.00613 & 0.00047 & 0.00057 & 0.00004 & 0.00006 \\
        Bias & -0.00323 & 0.00750 & 0.00221 & 0.00362 & 0.00053 & -0.00157 \\
        \bottomrule
    \end{tabular}
    \caption{Summary statistics of 100 estimates of $\gamma_1$ and $\gamma_2$ with different $K$.}
    \label{tab:beta_estimates}
\end{table}

\section{Confidence intervals using CLT for martingales} \label{S:CI}

Here we use the CLT for martingales in the evolution equation to obtain confidence limits for parameters.
Re-writing the martingale in \eqref{Evo1} and using an auxiliary result \cite[Proposition 26]{FanEtal20} of the CLT of the population process, we have 
\begin{equation}\label{CLT1}
\sqrt{K}\Big((f,\bar A^K_T) - (f,\bar A_0^K) - \int_0^T (L_{\bar A_s^K}f,\bar  A_s^K) ds\Big) \stackrel{d}{\approx}  N(0,(V^f_T)^2),
\end{equation}
where $(V^f_T)^2$ is given by \eqref{QVM}. 
This gives 
\begin{equation}\label{CLT2}
P \left( \Big| (f,\bar A^K_T) - (f,\bar A_0^K) - \int_0^T (f' - f h_{\bar A_s^K} + f(0) b_{\bar A_s^K},\bar A_s^K)ds\Big|\le \frac{c_\alpha V_T^f}{\sqrt K} \right) \approx 1-\alpha,
\end{equation}
where $c_\alpha = z_{\alpha/2}$ the upper percentage point of a standard normal distribution. 
More generally, for test functions of two variables, we have 
\begin{multline}\label{CLT3}
P\bigg(\Big|(f_T,\bar A^K_T) - (f_0,\bar A_0^K) - \int_0^T (\partial_xf_s + \partial_sf_s - f_sh_{\bar A^K_s} + f_s(0) b_{\bar A^K_s}, \bar A^K_s) ds \Big|  
\le \frac{ c_\alpha V_T^f }{ \sqrt K } \bigg) \\\approx 1-\alpha,
\end{multline}
where 
\begin{equation}\label{CLTV}
(V_T^f)^2 = \int_0^T (f^2_s(0) b_{\bar A_s} + h_{\bar A_s}f_s^2, \bar A_s)ds.
\end{equation}

Note that $V^f_T$ also involves unknown parameters as well as the limiting process $\bar A$.  Replacing $\bar A$ with its pre-limit $\bar A^K$ yields inequalities for the approximate confidence regions given below.
\begin{eqnarray}\label{CR}
\Big|(f_T,\bar A^K_T) - (f_0,\bar A_0^K) - \int_0^T (\partial_xf_s + \partial_sf_s - f_sh_{\bar A^K_s} + f_s(0) b_{\bar A^K_s}, \bar A^K_s) ds \Big| &\le &\nonumber
\\
\frac{ c_\alpha }{ \sqrt K } \int_0^T (f^2_s(0) b_{\bar A_s} + h_{\bar A_s}f_s^2, \bar A^K_s)ds.&&
\end{eqnarray}

Another approximation is obtained by  completely replacing  $V^f_T$ by its estimator
\begin{equation}\label{V_hat}
\hat{V}_T^f = \int_0^T (\hat b_{T}f^2_s(0)  + \hat h_{T}f_s^2, \bar A^{K}_s)ds,
\end{equation}
where $\hat b_T$ and $\hat h_T$ here denote the estimates of $b$ and $h$.  Note that by continuity, $\hat{V}_T^f $ is   consistent estimator of $V_T^f $, as $K\to\infty$; so that for large $K$ there is little difference between the exact value  $V_T^f $ and its estimator $\hat{V}_T^f $. Therefore, this approach may be more suitable in practice, as it leads to simpler confidence regions. Naturally, this results in some loss of accuracy, which can be assessed in specific examples.

It is important to note that by replacing a parameter by its estimate in the inequality for the confidence region changes the original probability of that region. 

Essentially, the construction of confidence regions by using test functions is akin to that for the mean of multivariate normal distribution with unknown covariance matrix. 
We do not attempt to give a complete solution here, but merely suggest a practical  way to implement approximations.

To obtain  more  precise confidence regions for parameters, in subsequent research, we shall   use  as many as we need test functions, noting that for
a pair of test functions $f$ and $g$  the 
covariances  are given by the predictable quadratic covariation  (sharp bracket) formula  for martingales $M_t^f$ and $M_t^g$ in \cite{FanEtal20}
$$\langle M^f,M^g \rangle_t=\int_0^t (f(0)g(0)b_{\bar A_s}+h_{\bar A_s}fg, \bar A_s)ds.$$

\subsection{Constant parameters} \label{CLT: const}

For the case of constant parameters, recall from Section \ref{S:const} the estimator of $h$ and $b$: 
\begin{align*}
\hat h_T &= \frac{(x,\bar{A}_0^K) - (x,\bar{A}_T^K) + \int_0^T (1,\bar{A}^{K}_s) ds}{\int_0^T (x,\bar{A}^{K}_s) ds}, \\
\hat b_T &= \frac{ (1,\bar{A}_T^K) - (1,\bar{A}_0^K) + \hat h \int_0^T (1,\bar{A}^{K}_s) ds}{\int_0^T (1,\bar{A}^{K}_s) ds};
\end{align*}
and
$$(V_T^f)^2 = \int_0^T (bf^2_s(0) + hf_s^2, \bar A_s)ds.$$

Taking $f(x)=x$ in \eqref{CLT2}  eliminates $b$. 
We obtain the following quadratic inequality for $h$  
\begin{equation}\label{Inh}
\Big|(x,\bar{A}^K_T) - (x,\bar{A}_0^K) - \int_0^T (1,\bar{A}^{K}_s)ds + h \int_0^T (x,\bar{A}^{K}_s) ds  \Big| \\
\le \frac{c_\alpha}{\sqrt K}  \sqrt{h} \sqrt{\int_0^T (x^2,\bar{A}^{K}_s) ds}.
\end{equation}
Solving it gives a confidence interval for $h$: 
\begin{equation}\label{CI_h_direct}
\left(\hat h_T + \frac{c^2_\alpha\int_0^T (x^2,\bar{A}^{K}_s) ds}{2K (\int_0^T (x,\bar{A}^{K}_s) ds)^2}\right) \pm  \frac{c_\alpha \sqrt{\int_0^T (x^2,\bar{A}^{K}_s) ds}}{\sqrt{K} \int_0^T (x,\bar{A}^{K}_s) ds}\sqrt{\hat h_T +   \frac{c^2_\alpha\int_0^T (x^2,\bar{A}^{K}_s) ds}{4K (\int_0^T (x,\bar{A}^{K}_s) ds)^2}}.
\end{equation}
Next, take $f(x)=1$ in \eqref{CLT2}, and note that $(V_T^1)^2 = (b+h)\int_0^T (1, \bar A_s)ds$.
We obtain
\begin{equation}\label{bh}
\Big| (1,\bar A^K_T) - (1,\bar A_0^K) + (b-h)\int_0^T (1,\bar A_s^K)ds\Big|\le \frac{c_\alpha}{\sqrt K}\sqrt {b+h}\sqrt{\int_0^T (1, \bar A^K_s)ds}.
\end{equation}
Inequalities \eqref{Inh} and \eqref{bh} define the confidence region for $h$ and $b$.

A naive approximate confidence interval for $b$ can be constructed by replacing $h$ by its estimator in the above inequality \eqref{bh}, reducing it to one-dimensional inequality.
Using estimate $\hat h_T$, we can solve the following inequality for $b$
$$\sqrt K\Big|(1,\bar{A}^K_T)  - (1,\bar{A}^K_0)  - (b - \hat h_T) \int_0^T (1,\bar{A}^{K}_s) ds  \Big| \\
\le c_\alpha  \sqrt{b+\hat h_T} \sqrt{\int_0^T (1,\bar{A}^{K}_s) ds} $$
and obtain a confidence interval of $b$: 
\begin{equation}\label{CI_b_direct}
 \left(\hat b_T + \frac{c^2_\alpha}{2K \int_0^T (1,\bar{A}^{K}_s) ds }\right) \pm  \frac{c_\alpha }{\sqrt{K\int_0^T (1, \bar{A}^{K}_s) ds} }\sqrt{\hat b_T + \hat h_T +   \frac{c^2_\alpha}{4K \int_0^T (1,\bar{A}^{K}_s) ds}}.
\end{equation}
Note that  confidence intervals in \eqref{CI_h_direct} and \eqref{CI_b_direct} are 
asymptotically centered at the estimators $\hat h_T$ and $\hat b_T$, with a vanishing shift of order $1/K$. 

Second approach is when we replace  $V^f_T$ by its estimate $\hat V^f_T$, given by \eqref{V_hat} .
Taking $f(x)=x$ in \eqref{CLT2}  we obtain a confidence interval of $h$,
\begin{equation*} 
\hat h_T \pm c_\alpha \frac{\sqrt{\hat h_T} \sqrt{\int_0^T (x^2,\bar{A}^{K}_s) ds}}{\sqrt{K}\int_0^T (x,\bar{A}^{K}_s) ds}.
\end{equation*}
Taking $f(x)=1$ we obtain a confidence interval of $b$,
\begin{equation*}
\hat b_T \pm c_\alpha \frac{ \sqrt{\hat b_T + \hat h_T} }{ \sqrt{K} \sqrt{\int_0^T (1,\bar{A}^{K}_s) ds}}.
\end{equation*}

\begin{figure}[h]
    \centering
    \begin{minipage}{0.43\textwidth}
        \centering
        \includegraphics[width=6cm]{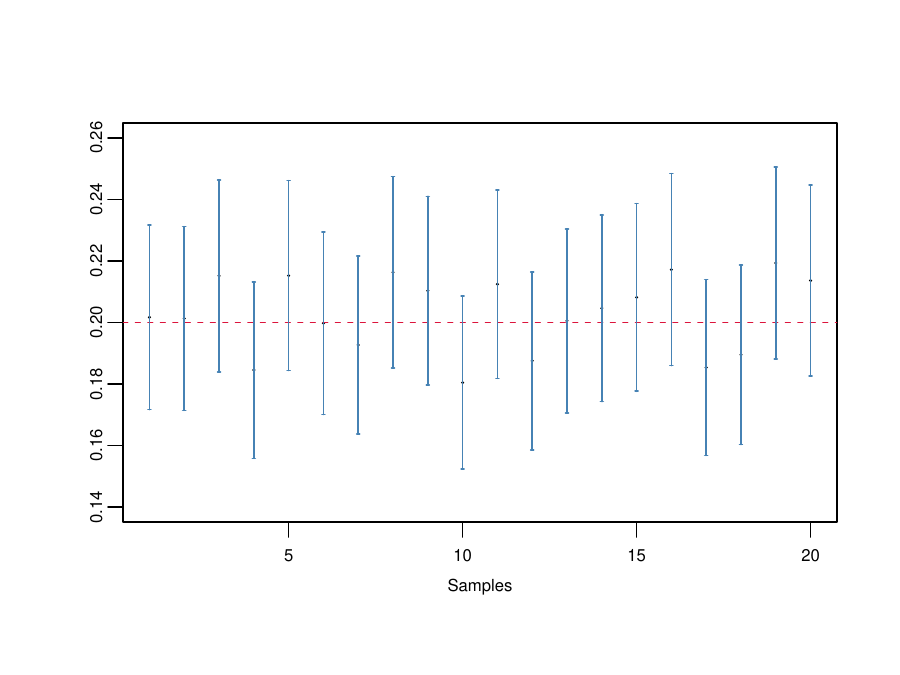} 
    \end{minipage}
    \hfill
    \begin{minipage}{0.5\textwidth} 
        \centering
        \includegraphics[width=6cm]{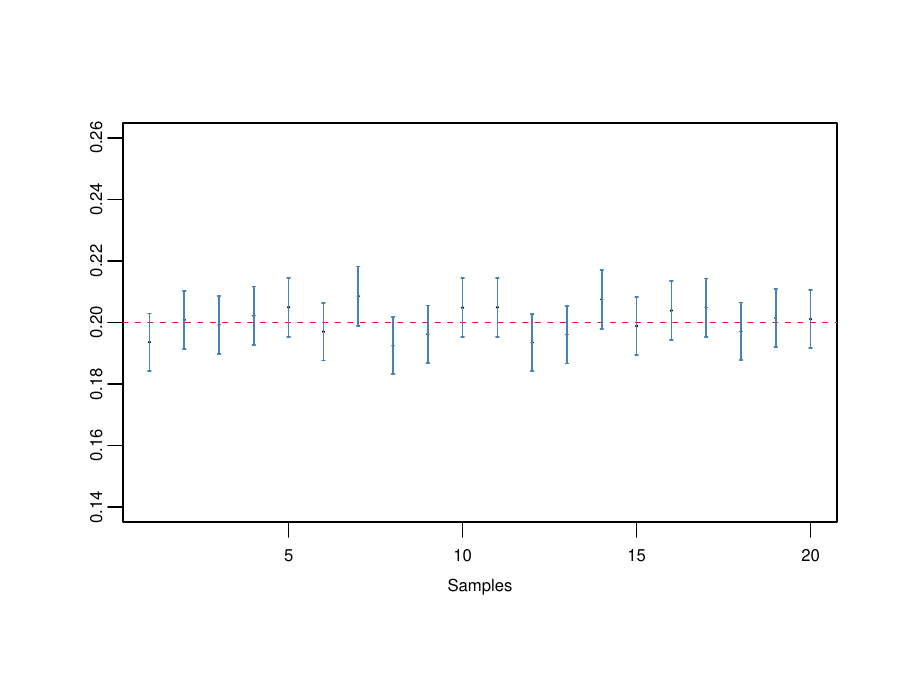} 
    \end{minipage}
    \begin{minipage}{0.43\textwidth} 
        \centering
        \includegraphics[width=6cm]{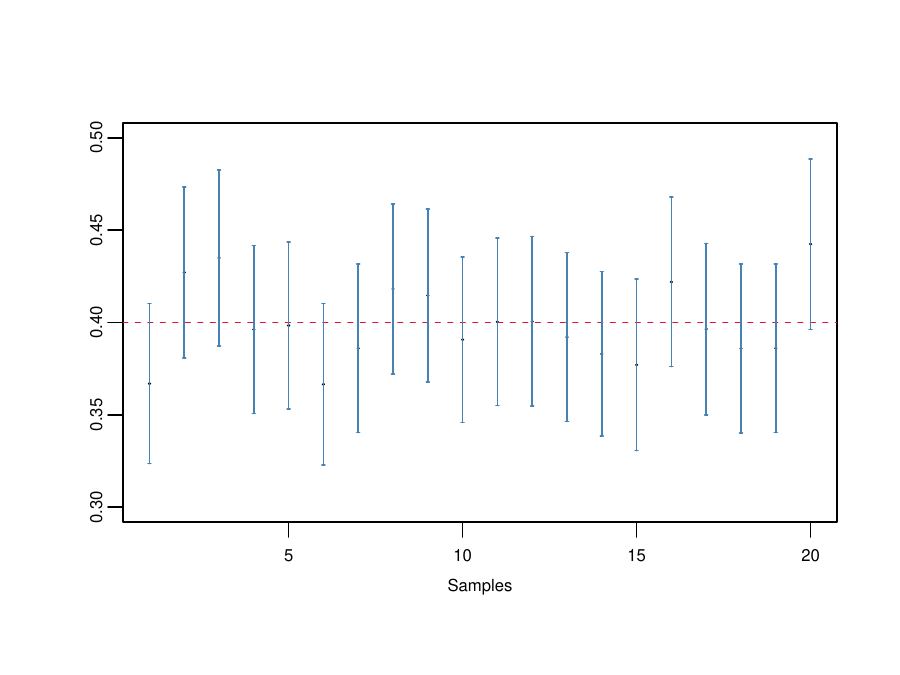} 
    \end{minipage}
    \hfill
    \begin{minipage}{0.5\textwidth}
        \centering
        \includegraphics[width=6cm]{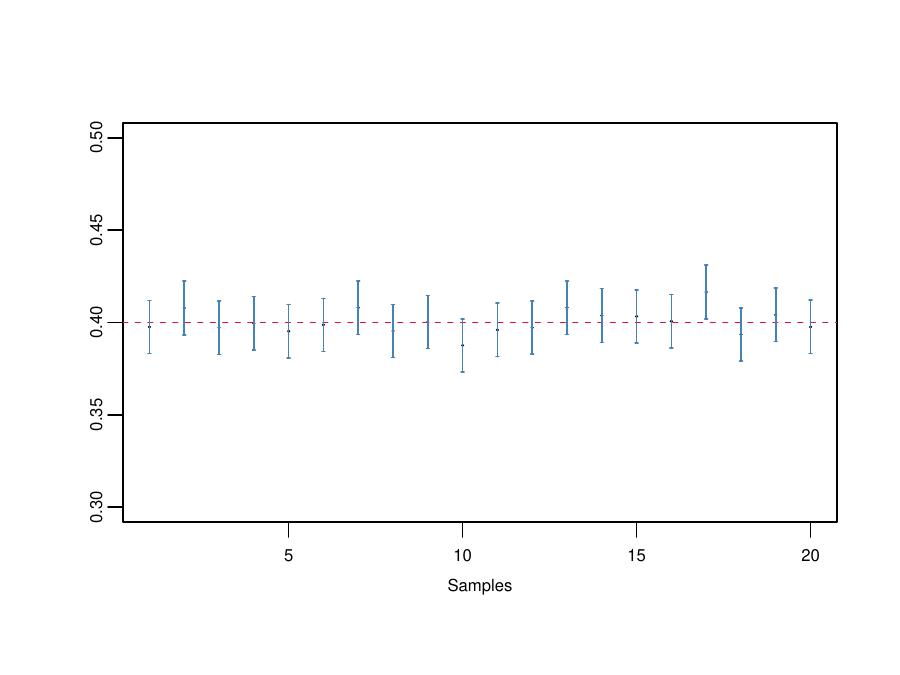} 
    \end{minipage}
    \caption{Confidence intervals of $h$ (top row) and $b$ (bottom row) in 20 samples with $K=1000$ (left) and $K=10000$ (right).}
    \label{fig:CI_const}
\end{figure}

Figure \ref{fig:CI_const} shows the confidence intervals of $h$ and $b$ 
for different $K$ values using the direct approach, \eqref{CI_h_direct} and \eqref{CI_b_direct}.  
These were obtained based on the same parameter values as in the numerical examples in Section \ref{S:const}.
As expected, shorter intervals are realised for larger $K$.

\subsubsection{Comparison with Classical result}

Our estimation produces classical result  for constant rates. We can see this as follows.  

Consider a pure birth process with $b>0$ and $h=0$. 
Taking $f=1$, \eqref{EvoA} gives 
$$(1,\bar A_t) = (1,\bar A_0) + b\int_{0}^{t} (1,\bar{A}_s)ds,$$
which results in estimator 
$$\hat b_T = \frac{(1,\bar A^{K}_T)-(1,\bar A^{K}_0)}{\int_{0}^{T} (1,\bar A^{K}_s)ds}.$$
This estimator coincides with the Maximum Likelihood Estimator in \cite{Keiding74}. 

Moreover, from \eqref{CLT2} and replacing $\bar A$ with $\bar A^K$ in $V_T^1$,
we have
\begin{align*}
1-\alpha &\approx 
P\left(\sqrt K\Big|(1,\bar A^K_T)-(1,\bar A_0^K)-b\int_0^T (1, \bar A^{K}_s)ds\Big|\le c_\alpha \sqrt{ b \int_0^T (1, \bar A^{K}_s)ds} \right) \\
&=P\left(\sqrt{\frac{K\int_{0}^{T} (1,\bar A^{K}_s)ds}{b}} \; \Big|\hat b_T - b\Big|\le  c_\alpha  \right). 
\end{align*}
Thus, 
$$\sqrt{\frac{K \int_{0}^{T} (1,\bar A^{K}_s)ds}{b}} \left(\hat b_T - b\right) \stackrel{d}{\approx} N(0, 1),$$
for any $T$, 
which is consistent with   \cite[Theorem 3.5(a)]{Keiding74}.

\subsection{Parameters depend only on population}

Suppose $b_A = \eta (1_{J_1},A)$ and $h_A = \lambda (1_{J_2},A)$ as in Section \ref{S:PopDep}. 
Recall that the estimators of $\lambda$ and $\eta$ are
\begin{equation*}
\hat \lambda_T = \frac{ (x, \bar A^K_0) - (x, \bar A^K_T) + \int_0^{T} (1,\bar A^K_s) ds}{\int_{0}^{T} (x,\bar A^K_s) (1_{J_2},\bar A^K_s) ds},
\end{equation*}
and
\begin{equation*}
\hat \eta_T = \frac{(1,\bar A^K_T) - (1, \bar A^K_0) + \hat \lambda_T \int_{0}^{T} (1_{J_2},\bar A^{K}_s)(1,\bar A^K_s)ds}{\int_{0}^{T}(1,\bar A^K_s)(1_{J_1}, \bar A^K_s)ds}.
\end{equation*} 

Taking $f(x,t)=x$, 
$$(V^{x}_T)^2  =  \lambda \int_0^T  (1_{J_2}, \bar A_s)(x^2,\bar{A}_s) ds. $$
Replacing $\bar A$ with $\bar A^K$ in $V^{x}_T$, from \eqref{CLT2} a confidence interval of $\lambda$ is obtained by solving 
\begin{equation*}
\sqrt K\Big| \lambda - \hat \lambda_T \Big| \\
\le c_\alpha  \frac{\sqrt{\lambda} \sqrt{\int_0^T (1_{J_2}, \bar A^K_s) (x^2,\bar{A}^{K}_s) ds}}{\int_{0}^{T}(1_{J_2},\bar A^K_s) (x,\bar A^K_s)  ds}.
\end{equation*}
This gives 
\begin{equation}\label{CI_popdep_h_direct}
\left(\hat \lambda_T + \frac{c^2_\alpha I_{J_2}^{x^2}}{2K \big( I_{J_2}^{x}\big)^2}\right) \pm  \frac{c_\alpha \sqrt{I_{J_2}^{x^2}}}{\sqrt{K} I_{J_2}^{x}}\sqrt{\hat \lambda_T +   \frac{c^2_\alpha I_{J_2}^{x^2} }{4K \big(I_{J_2}^{x}\big)^2}},
\end{equation}
where 
$$I^f_J := I^f_J(T) =  \int_0^T  (1_{J}, \bar A^K_s)(f,\bar{A}^{K}_s) ds.$$
Similarly,  a confidence interval of $\eta$ is obtained by taking $f(x)=1$, which yields   
$$\left(\hat \eta_T + \frac{c^2_\alpha }{2K I_{J_1}^1}\right) \pm  \frac{c_\alpha }{ \sqrt{K } I_{J_1}^1}\sqrt{\hat \eta_T I_{J_1}^1+ \hat \lambda_T I_{J_2}^1 + \frac{c^2_{\alpha}}{4K}} .$$

Alternatively, we can replace the unknown parameters in $V^f_T$ with their estimates. 
Then, taking $f(x)=x$, we get a confidence interval of $\lambda$:
\begin{equation}\label{CI_popdep_h_approx}
\hat \lambda_T \pm c_\alpha \frac{\sqrt{\hat \lambda_T \int_0^T  (1_{J_2}, \bar A^K_s)(x^2,\bar{A}^{K}_s) ds}}{\sqrt{K}\int_{0}^{T} (1_{J_2},\bar A^K_s) (x,\bar A^K_s) ds};
\end{equation}
and taking $f(x)=1$ gives a confidence interval of $\eta$:
$$\hat \eta_T \pm c_\alpha \frac{\sqrt{\hat \eta_T  \int_0^T (1_{J_1}, \bar A^K_s)(1,\bar{A}^{K}_s) ds + \hat \lambda_T \int_0^T (1_{J_2}, \bar A^K_s)(1,\bar{A}^{K}_s) ds} }{\sqrt{K}\int_{0}^{T} (1,\bar A^K_s) (1_{J_1},\bar A^K_s) ds}.$$

\begin{figure}[h]
    \centering
    \begin{minipage}{0.43\textwidth}
        \centering
        \includegraphics[width=6cm]{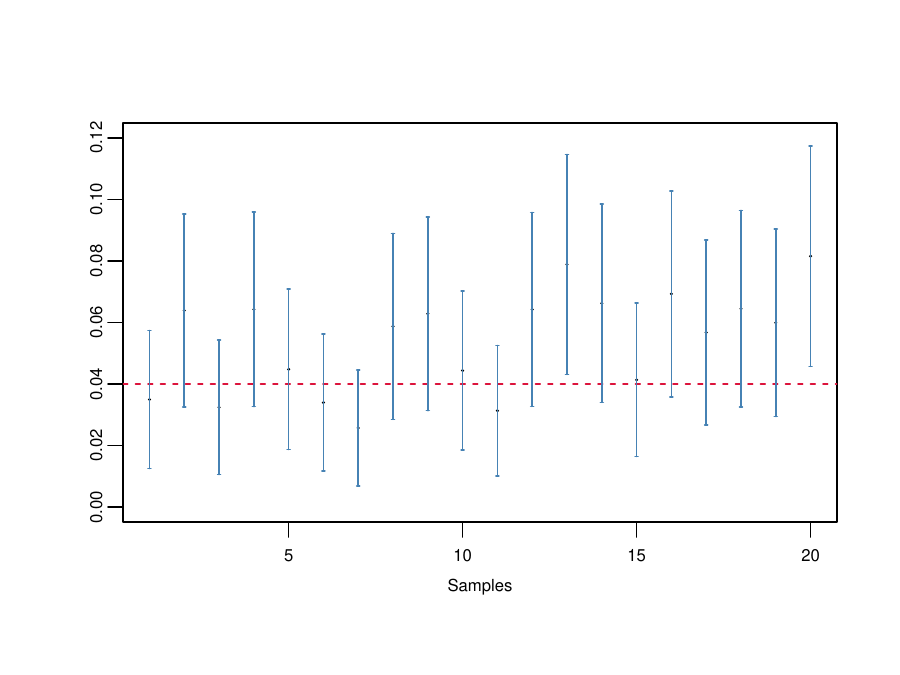} 
    \end{minipage}
    \hfill 
    \begin{minipage}{0.5\textwidth}
        \centering
        \includegraphics[width=6cm]{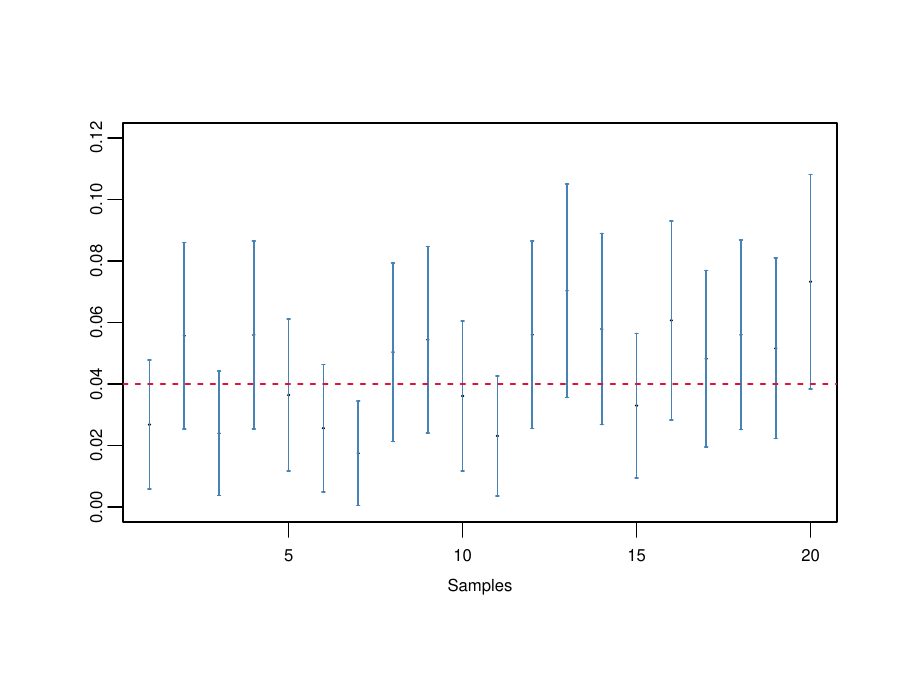} 
    \end{minipage}
    \caption{Confidence intervals of $\lambda$ in 20 samples with $K=1000$ using the direct approach eq. \eqref{CI_popdep_h_direct} (left) and the approximate approach eq. \eqref{CI_popdep_h_approx} (right).}
    \label{fig:CI_PopDep_h}
\end{figure}

Figure \ref{fig:CI_PopDep_h} shows confidence intervals of $\lambda$ obtained using the two approaches from the same sample. 
These were obtained based on the same parameter values as in the numerical examples in Section \ref{S:PopDep}.
Note that the direct approach resulted in higher intervals. 

\subsection{Parameters depend only on age} \label{S:CLT_AgeDep}

Suppose
\begin{equation*}
h(x)=\sum_{i=1}^n h_i1_{B_i}(x) \quad \textnormal{and} \quad b(x)=\sum_{i=1}^n b_i1_{B_i}(x) 
\end{equation*}
as in Section \ref{sec: age_depend}. 
In this case, we will  need test functions of two variables. 
Obtaining confidence intervals of $h_i$'s and $b_i$'s involves solving a system of inequalities. 

We provide a brief insight into the problem by considering the case $n=2$.
For $h_i$'s, take  $f(x) = x$ and $f_t(x)=xt$. We have
\begin{equation*}
 \left |\sum_{i=1}^2 h_i \int_0^T (x1_{B_i}(x),\bar A^K_s) ds +  (x,\bar A^K_T) - (x, \bar A^K_0) - \int_0^T (1,\bar A^K_s)ds \right | \le c_\alpha V_T^x/\sqrt K,
\end{equation*}
and
\begin{equation*}
\left| \sum_{i=1}^2 h_i \int_0^T s(x1_{B_i}(x),\bar A^K_s) ds + T(x,\bar A^K_T) - \int_0^T s(1,\bar A^K_s) ds 
- \int_0^T  (x,\bar A^K_s) ds \right | \le c_\alpha V_T^{xt}/\sqrt K .
\end{equation*}

The direct approach with $\bar A^K$ in $V^f_T$ gives a system of nonlinear inequalities. 
A confidence region for $\bm{h} = (h_1, h_2)$ is determined by identifying the feasible region of the above system of nonlinear inequalities.
Each inequality alone above  forms an elliptical region in some space. This happens when the constraints define an ellipse (in 2D) as a feasible region. 

The confidence region of $\bm{b} = (b_1, b_2)$ can be obtained in a similar way by using  estimates of $h_i$'s, and taking $f(x,t)=1$ and $f(x,t)=t$.

Alternatively, using the estimate $\hat V^f_T$ given in \eqref{V_hat} gives a system of linear inequalities.

\begin{figure}[H]
    \centering
    \begin{minipage}{0.43\textwidth} 
        \centering
        \includegraphics[width=6cm]{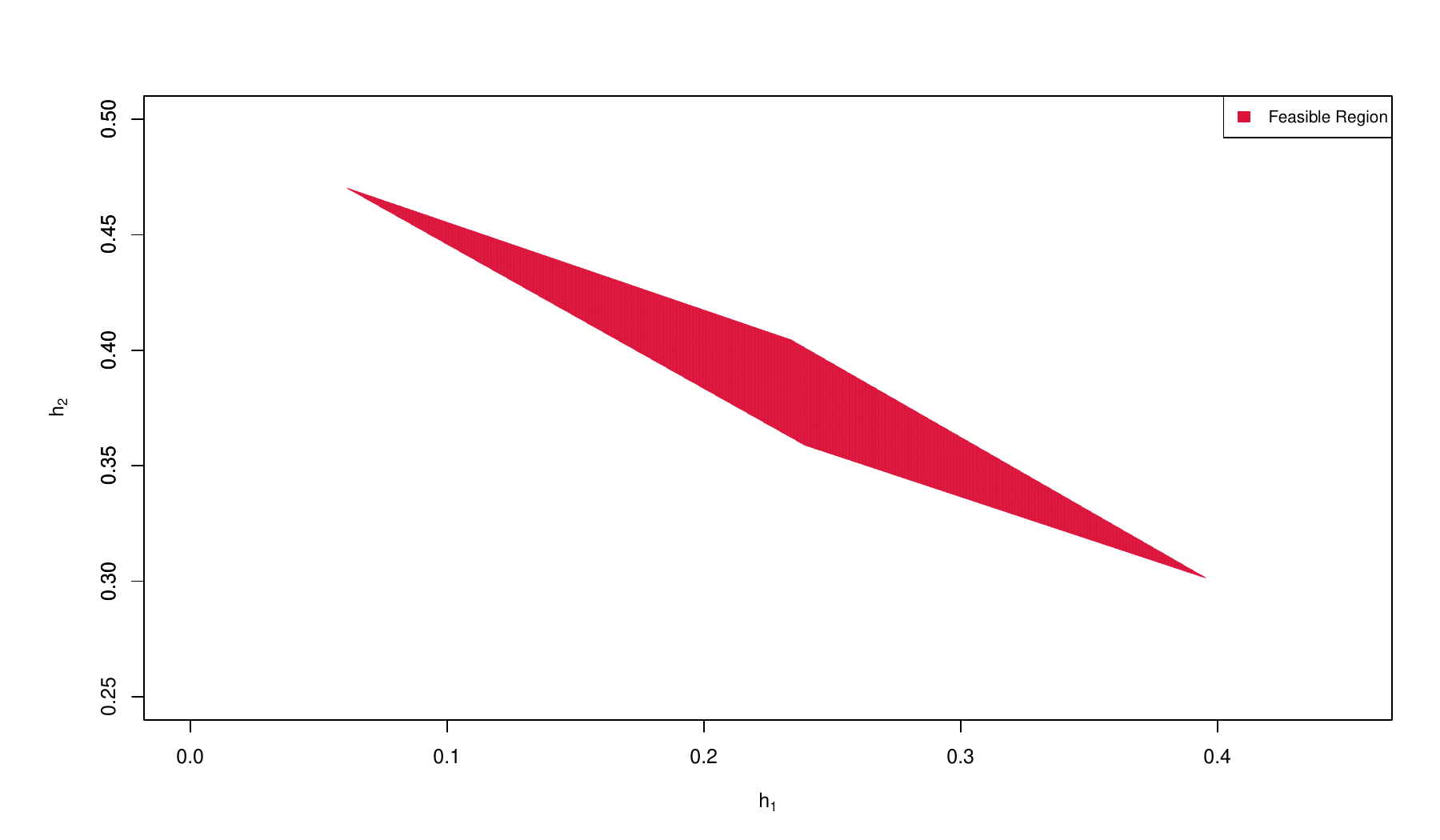} 
    \end{minipage}
    \hfill
    \begin{minipage}{0.5\textwidth} 
        \centering
        \includegraphics[width=6cm]{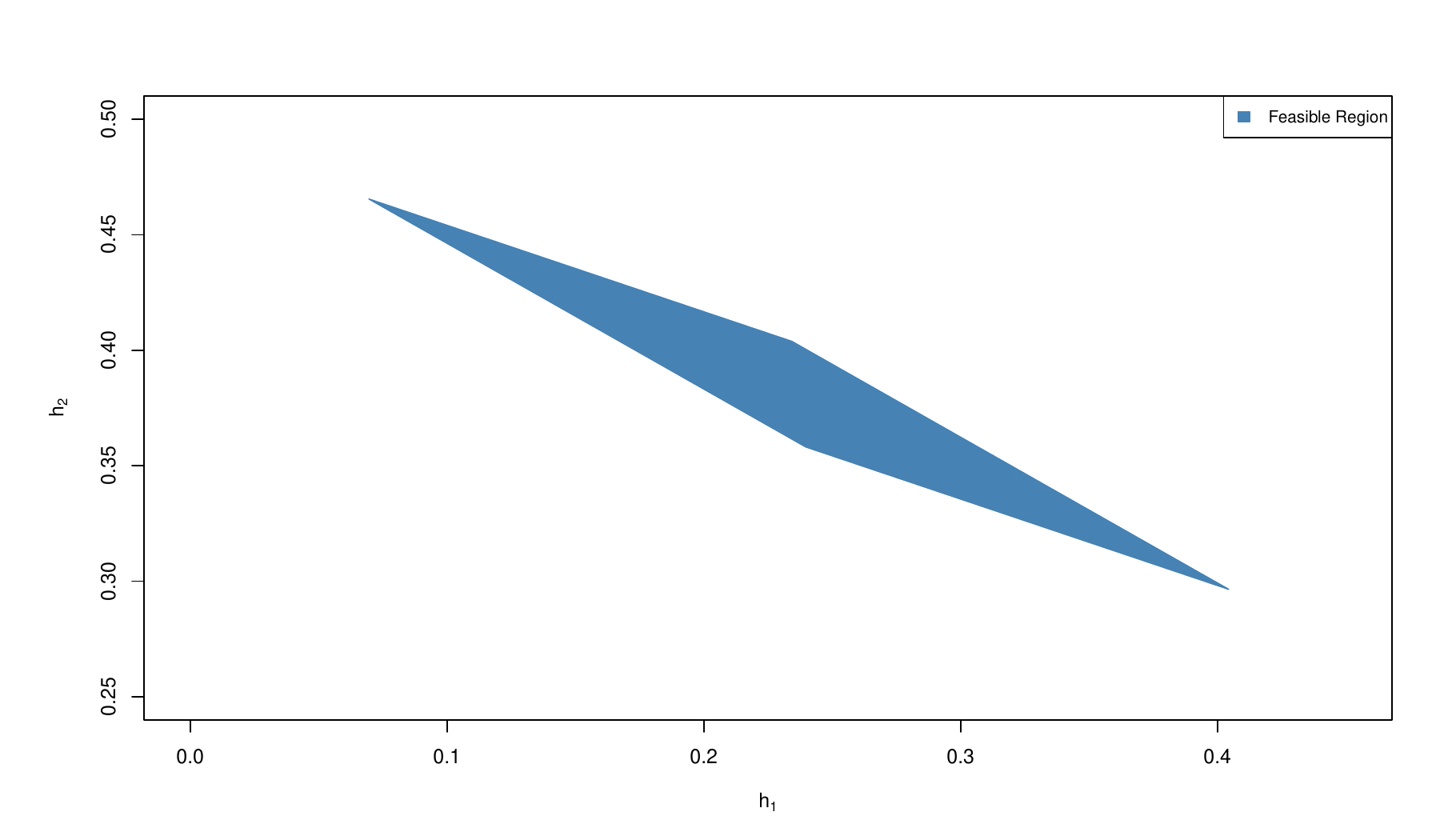} 
    \end{minipage}
    \begin{minipage}{0.45\textwidth}
        \centering
        \includegraphics[width=6cm]{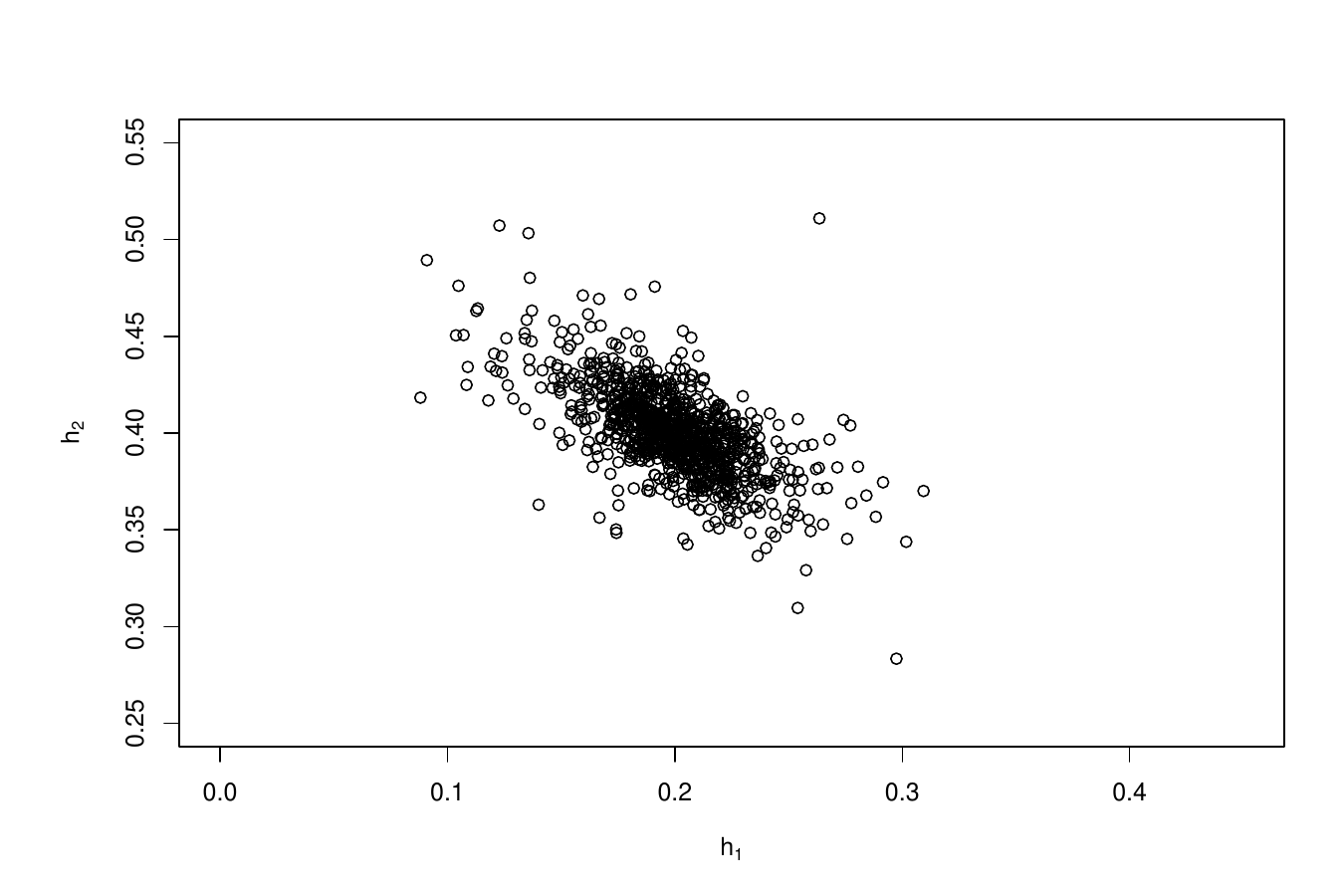} 
    \end{minipage}
    \caption{Confidence regions of $\bm{h}$ in one sample with $K=10000$ using the direct approach (top left) and the approximate approach (top right). Point estimates of $(h_1, h_2)$ in 100 samples for $K=10000$ (bottom).}
    \label{fig:CI_AgeDep_h}
\end{figure}

Figure \ref{fig:CI_AgeDep_h} shows confidence region of $(h_1,h_2)$ obtained using the two approaches from the same sample. These were obtained based on the same parameter values as in the numerical examples in Section \ref{sec: age_depend}.
As a comparison, a plot of 100 point estimates of $(h_1, h_2)$ in 100 samples for $K=10000$ is also given.

\subsection{Parameters depend on population and age}

The general case where  parameters $h$ and $b$ depend on both $A$ and $x$ can be dealt with in a similar way as in Section \ref{S:CLT_AgeDep}. 
In particular, when $b_A$ and $h_A$ take the forms of \eqref{eq: popage_h} and \eqref{eq: popage_b} as in Section \ref{sec: popage}: 
\begin{gather*}
h_A(x)=\alpha_1 (1_{J},A) 1_{B_1}(x) + \alpha_2 (1_{J},A) 1_{B_2}(x),\\
b_A(x)=\gamma_1 (1_{J},A )1_{B_1}(x) + \gamma_2 (1_{J},A) 1_{B_2}(x). 
\end{gather*}
This can be generalised using the same idea. 

Taking  $f_t(x) = x$ and $f_t(x) = xt$ in \eqref{CLT3}, we have a system of inequalities: 
\begin{gather*}
\makebox[\textwidth][c]{$\displaystyle\left | \sum_{i=1}^2 \alpha_i \int_0^T (1_J,\bar A^K_s) (x1_{B_i}(x),\bar A^K_s) ds + (x,\bar A^K_T) - (x,\bar A_0^K)  - \int_0^T (1,\bar A^K_s)ds \right | \le c_\alpha \frac{V_T^x}{\sqrt K},$} \\
\makebox[\textwidth][l]{$\displaystyle\left| \sum_{i=1}^2 \alpha_i \int_0^T s (1_J,\bar A_s^K) (x1_{B_i}(x),\bar A^K_s) ds + T(x,\bar A^K_T) - \int_0^T s(1,\bar A^K_s) ds 
- \int_0^T  (x,\bar A^K_s) ds \right |$} \\
\makebox[\textwidth][r]{$\displaystyle\le c_\alpha \frac{V_T^{xt}}{\sqrt K} .$}
\end{gather*}
With $\bar A^K$ in $V^f_T$, solving this system of nonlinear equations, we obtain a confidence region of $\bm\alpha= (\alpha_1, \alpha_2)$.
The same with  $f_t(x) = 1$ and $f_t(x) = t$ gives a confidence region of $\bm\gamma= (\gamma_1, \gamma_2)$.

Alternatively, using $\hat V_f^T$ in \eqref{V_hat}, the confidence regions of $\bm\alpha$ and $\bm\gamma$ 
can be obtained through a system of linear inequalities.

%
\begin{acknowledgement} We thank the referee for the  comments.\\
This research was supported by the
Australian Research Council grant DP220100973.
\end{acknowledgement}

\bigskip

\end{document}